\documentclass{amsart}
\usepackage{graphicx}
\newtheorem{theorem}{Theorem}[section]
\newtheorem{corollary}[theorem]{Corollary}

\theoremstyle{definition}

\theoremstyle{remark}
\newtheorem{rem}[theorem]{Remark}
\numberwithin{equation}{section}

\begin{document}

\title[]{A note on  recurrence sequences }%
\author{S. Amghibech}%
\address{}%
\email{amghibech@hotmail.com}%

\thanks{}%
\subjclass{}%
\keywords{}%

\begin{abstract}
In this note we provide a simple  formula of general term of
recurrent sequence.
\end{abstract}
\maketitle

Suppose that the sequence $\{u_{h}\}$ satisfies the relation
\begin{equation*}\label{recurrence}
u_{h+n}=\sum_{i=1}^{n}s_{i}u_{n+h-i} \ \ \ \ \text{for}\ \ \
h=0,1,2,\ldots
\end{equation*}

for some complex numbers $s_1,\ldots s_n$, $s_{n}\not=0$ and for
$h=0,1,\ldots $. Taking $h=0$ we see that $u_{n}$ is in the ring
$\mathbb{Z}[u_{0},\ldots u_{n-1},s_{1},\ldots s_{n}]$, by using
the induction argument we show that all the terms in the sequence
belong to this ring.

Let
$$
X^{n}-\sum_{i=1}^{n}s_{i}X^{n-i}=\prod_{i=1}^{m}(1-\alpha_{i}X)^{n_i}
$$
with distinct nonzero roots $\alpha_{i}$. It well known that there
exist $A_{i}(X)$ which are polynomials of degree $n_{i}-1$ for
positive integers $n_{i}$,  $i\in [1,m]$ such that
$$
u_{h}=\sum_{i=1}^{m}A_{i}(h)\alpha_{i}^{h}
$$
 for $h = 0, 1, \ldots$, see for example \cite{Mye}, for an original work
  on recurrent sequences see \cite{Mon}. More
 concretely, in this note we give a short proof of this result and the
  explicit formula $\{u_{h}\}$.

Consider the vector space $E$ of all sequences $\{v_{h}\}$ such
that
$$
v_{h+n}=\sum_{i=1}^{n}s_{i}v_{n+h-i}\ \ \text{for}\ \
h=0,1,\ldots.
$$
Remark that $\dim E=n$ and that $e^{j,i}=\{h(h-1)\ldots
(h+1-j)\alpha_i^{h-j}\}$ for $j=0,\ldots ,n_i-1$, $i=1,\ldots, m$
is a basis of $E$.

Define the projection $\Psi_{k} : E\longrightarrow
\mathbb{C}^{n+1}$ for $k\geq n$ by
$$
\Psi_{k}(\{v_{h}\}):= (v_{0},\ldots ,v_{n-1}, v_{k} )^{t}.
$$
From the fact that $ \dim E=n$ we get the vectors
$\Psi_{k}(\{u_{h}\})$, $\Psi_{k}(e^{j,i})$, $j=0,\ldots ,n_i-1$,
$i=1,\ldots, m$ are linearly dependent, thus
$$
\det [\Psi_{k}(\{u_{h}\}),\Psi_{k}(e^{j,i})_{0\leq j\leq n_i-1,
1\leq i\leq m} ]=0
$$
which gives
$$
(-1)^{n}u_{h}\det[\Psi_{n}(e^{j,i})_{0\leq j\leq n_i-1, 1\leq
i\leq m} ]+ \det[f, \Psi_{k}(e^{j,i})_{0\leq j\leq n_i-1, 1\leq
i\leq m} ]=0
$$
with $f=(u_{0},\ldots ,u_{n-1}, 0 )^{t}$, from the fact that
$e^{j,i}=\{h(h-1)\ldots (h+1-j)\alpha_i^{h-j}\}$ for $j=0,\ldots
,n_i-1$, $i=1,\ldots, m$ is a basis of $E$ we get
$\det[\Psi_{n}(e^{j,i})_{0\leq j\leq n_i-1, 1\leq i\leq m}
]\not=0$ thus
$$
u_{h}=(-1)^{n+1} \frac{\det[f, \Psi_{k}(e^{j,i})_{0\leq j\leq
n_i-1, 1\leq i\leq m} ]}{\det[\Psi_{n}(e^{j,i})_{0\leq j\leq
n_i-1, 1\leq i\leq m} ]}.
$$
Also  the (ordinary) generating function of the sequence
$\{u_{h}\}$, is given by
$$
\sum_{h=0}^{\infty}u_{h}x^{h} =(-1)^{n+1} \frac{\det[f,
f_{j,i}(x)_{0\leq j\leq n_i-1, 1\leq i\leq m} ]}
{\det[\Psi_{n}(e^{j,i})_{0\leq j\leq n_i-1, 1\leq i\leq m} ]}
$$
with
$$
f_{j,i}(x)=(\Psi_{n}(e^{j,i})^{t},\frac{j!x^{j}}{(1-\alpha_ix)^{j+1}})^{t}
$$
for $0\leq j\leq n_i-1$ and $1\leq i\leq m$. Of course the vectors
are ordered in the same way in the two determinants.

\vskip 5mm

Let $I$ be finite subset of $\mathbb{N}$ and define the projection
$\Psi_{I} : E\longrightarrow \mathbb{C}^{I}$ by
$$
\Psi_{I}(\{v_{h}\}):= (v_{l})_{l\in I} ^{t}.
$$
By using a similar argument as above we get the following theorem
\begin{theorem} Let $I$ be a subset of $n$ elements of $\mathbb{N}$, let $h\not\in
I$, put $J=I\cup \{h\}$. Then we have
$$
u_{h}=(-1)^{n+1} \frac{\det[((\Psi_{I}(\{u_{h}\}))^{t},0)^{t},
\Psi_{J}(e^{j,i})_{0\leq j\leq n_i-1, 1\leq i\leq m}
]}{\det[\Psi_{I}(e^{j,i})_{0\leq j\leq n_i-1, 1\leq i\leq m} ]}.
$$
if $\det[\Psi_{I}(e^{j,i})_{0\leq j\leq n_i-1, 1\leq i\leq m}
]\not=0$.
\end{theorem}
\vskip 5mm

 In the case $m=n$ we have $e^{0,i}=\{\alpha_i^{h}\}$ is
a basis of $E$ then we get following corollary
\begin{corollary}\label{recurrence}
For all $h$ we have
$$
u_{h}=(-1)^{n+1}\prod_{1\leq i<j\leq
n}(\alpha_{j}-\alpha_{i})^{-1}
\begin{vmatrix}
  u_{0} & 1 & \ldots & 1 \\
  u_{1} & \alpha_{1} & \ldots & \alpha_{n}\\
 \vdots & \vdots & \ddots & \vdots \\
  u_{n-1}& \alpha_{1}^{n-1}& \ldots & \alpha_{n}^{n-1} \\
  0     & \alpha_{1}^{h} & \ldots & \alpha_{n}^{h} \\
\end{vmatrix}
$$
\end{corollary}
Proof. Follows from  the theorem by remarking that
$I=\{0,1,\ldots,n-1\}$, and $\det[\Psi_{I}(e^{0,i})_{ 1\leq i\leq
m} ]$ is the Van der Mond determinant with $m=n$, when $h\not \in
I$. And the equality is also true when $h\in I$, which achieve the
proof.

 \vskip 5mm

In the next corollary  gives $u_h$ as a function of an $n$ fixed
terms of the sequence $\{u_h\}$
\begin{corollary}Let $k_{0}<k_{1}<\cdots <k_{n-1}$ a sequence of positive
integers. Then for all $h$ we have
$$
u_{h}=(-1)^{n+1}
\begin{vmatrix}
    \alpha_{1}^{k_{0}} & \ldots & \alpha_{n}^{k_{0}}\\
    \alpha_{1}^{k_{1}} & \ldots & \alpha_{n}^{k_{1}}\\
      \vdots           &  \ddots & \vdots \\
   \alpha_{1}^{k_{n-1}}& \ldots & \alpha_{n}^{k_{n-1}} \\
\end{vmatrix}^{-1}
\begin{vmatrix}
  u_{k_0} & \alpha_{1}^{k_{0}} & \ldots & \alpha_{n}^{k_{0}} \\
  u_{k_1} & \alpha_{1}^{k_{1}} & \ldots & \alpha_{n}^{k_{1}}\\
 \vdots & \vdots & \ddots & \vdots \\
  u_{k_{n-1}}& \alpha_{1}^{k_{n-1}}& \ldots & \alpha_{n}^{k_{n-1}} \\
  0     & \alpha_{1}^{h} & \ldots & \alpha_{n}^{h} \\
\end{vmatrix}
$$
when the expression has sense.
\end{corollary}

Proof.  Similar to the proof of Corollary \ref{recurrence}.
\begin{rem} $> rsolve( {f(0)=1, f(1)=1, f(n)=f(n-1)+f(n-2)},
f(n))$ gives the generale terms $f(n)$ in Maple.  
\end{rem}
\bibliographystyle{amsplain}
\bibliography{recu}
\end{document}